\documentclass[11pt,a4paper]{article}

\usepackage[a4paper,margin=2.45cm]{geometry}
\usepackage[T1]{fontenc}
\usepackage{lmodern}
\usepackage{amsmath,amssymb,amsthm,mathtools}
\usepackage{booktabs,longtable,array,multirow,tabularx}
\usepackage{float}
\usepackage{enumitem}
\usepackage{xcolor}
\usepackage{hyperref}
\usepackage{listings}
\usepackage{microtype}

\hypersetup{
  colorlinks=true,
  linkcolor=blue!55!black,
  citecolor=blue!55!black,
  urlcolor=blue!55!black,
  pdftitle={Small-LCM Distinct Covering Systems via Segmented Local Search},
  pdfauthor={Jiheng Zhang}
}
\setlist{nosep}
\renewcommand{\arraystretch}{1.15}
\newtheorem{theorem}{Theorem}
\newtheorem{conjecture}{Conjecture}

\theoremstyle{definition}
\newtheorem{definition}{Definition}

\newcommand{\Z}{\mathbb Z}
\newcommand{\D}{D}
\newcommand{\B}{\mathcal B}
\DeclareMathOperator{\lcm}{lcm}
\newcommand{\coverentry}[2]{\ensuremath{#1\,(#2)}}
\newcommand{\certhead}{\multicolumn{5}{c}{Entries are $a\,(d)$, meaning $x\equiv a\pmod d$}\\}
\lstdefinestyle{pseudo}{
  basicstyle=\small\ttfamily,
  columns=fullflexible,
  frame=single,
  breaklines=true,
  showstringspaces=false,
  numbers=left,
  numberstyle=\scriptsize\color{gray},
  keywordstyle=\bfseries\color{blue!55!black},
  commentstyle=\itshape\color{green!35!black},
  morekeywords={Input,Output,for,while,if,then,return,freeze,revisit}
}

\title{Small-LCM Distinct Covering Systems\\
via Segmented Local Search}

\author{
\begin{tabular}{c}
Jiheng Zhang\textsuperscript{1}\\[0.35em]
\textsuperscript{1}Institute of Mathematical Sciences, ShanghaiTech University
\end{tabular}
}
\date{November 1, 2026}

\begin{document}
\maketitle

\begin{abstract}
We construct distinct covering systems with minimum moduli $8,9,10,11$ and least common multiples $20160,75600,302400,3326400$, respectively, and conjecture that these values are minimal. For minimum modulus $8$, the exclusion computations leave only $10080$ and $15120$ unresolved below our construction. The constructions use segmented local search: moduli are introduced in increasing blocks, followed by joint reoptimization of selected blocks. We give complete residue lists and an independent integer-arithmetic verifier. Two computational studies, comprising $180$ initial runs and $130$ follow-up runs, compare initialization orders, revisit rules and seed-bank capacities under fixed scan limits.
\par\medskip
\noindent\textbf{Keywords:} covering systems; distinct moduli; least common multiple; local search; computational number theory
\end{abstract}

\section{Introduction}\label{sec:introduction}

Covering systems were introduced by Erd\H{o}s in 1950 \cite{Erdos1950}. For their history and subsequent developments, see \cite{PorubskySchonheim2002,Balister2024}. One basic question concerns the size of the minimum modulus in a distinct covering system. Hough proved that this modulus is bounded above by an absolute constant \cite{Hough2015}.

A related problem prescribes the minimum modulus and asks for the smallest possible least common multiple. For $m\ge2$, let $L_{\min}(m)$ denote the minimum LCM among distinct covering systems with minimum modulus $m$, with $L_{\min}(m)=\infty$ if no such system exists. 

Several small cases are settled. We have $L_{\min}(2)=12$; Dalton and Trifonov proved $L_{\min}(3)=120$ and $L_{\min}(4)=360$ \cite[Theorem~5]{DaltonTrifonov2022}. Klein established $L_{\min}(5)=1440$ and $L_{\min}(6)=5040$, and gave a covering with minimum modulus $7$ and least common multiple $15120$ \cite[Theorems~3--5]{Klein2026}. In previous joint work, my coauthor and I proved that $L_{\min}(7)=10080$ \cite{ZhangZhang2026}.

We give explicit upper bounds for the next four minimum moduli.

\begin{theorem}\label{prop:upper-bounds}
There exist distinct covering systems with minimum moduli $8,9,10,11$ and least common multiples $20160,75600,302400,3326400$, respectively. Consequently,
\[
  L_{\min}(8)\le20160,\quad L_{\min}(9)\le75600,\quad
  L_{\min}(10)\le302400,\quad L_{\min}(11)\le3326400.
\]
\end{theorem}
\begin{proof}
For each list in Appendix~\ref{app:constructions}, the verifier checks that the moduli are distinct, their minimum and LCM have the stated values, and every integer in $0,\ldots,L-1$ is covered. Periodicity then gives coverage of $\Z$.
\end{proof}

We conjecture that all four bounds are sharp.
\begin{conjecture}\label{conj:minimal-lcm}
\[
 L_{\min}(8)=20160,\qquad L_{\min}(9)=75600,\qquad
 L_{\min}(10)=302400,\qquad L_{\min}(11)=3326400.
\]
\end{conjecture}
Theorem~\ref{prop:upper-bounds} proves only the upper bounds. The conjecture requires excluding every smaller candidate LCM; unsuccessful heuristic searches do not provide such exclusions.

\begin{table}[htbp]
\centering
\caption{Our LCMs and earlier constructions. The last column gives the smallest earlier LCM located in the cited literature, excluding the author's work, as of November~1, 2026. It is not an exhaustive record claim; a dash denotes an unresolved comparison.}
\label{tab:headline-results}\label{tab:bosma}
\small
\begin{tabular}{@{}crrl@{}}
\toprule
$m$ & Our LCM & \shortstack{Bosma's greedy\\construction \cite{Bosma2006}} & \shortstack{Smallest prior LCM\\located (source)}\\
\midrule
8 & 20160 & 60480 & 50400 (Krukenberg; \cite[p.~3]{DaltonTrifonov2022})\\
9 & 75600 & 151200 & 151200 (Krukenberg; \cite[p.~3]{DaltonTrifonov2022})\\
10 & 302400 & 1663200 & 1663200 (Bosma \cite[p.~5]{Bosma2006})\\
11 & 3326400 & --- & ---\\
\bottomrule
\end{tabular}
\end{table}

Krukenberg's constructions have $L=2^5 3^2 5^2 7=50400$ for $m=8$ and $L=2^5 3^3 5^2 7=151200$ for $m=9$, as recorded by Dalton and Trifonov \cite[p.~3]{DaltonTrifonov2022}. Bosma's randomized greedy search gives the comparison periods in Table~\ref{tab:bosma}. Bickel et al.\ report $997920$ and $21621600$ for $m=10,11$ \cite[pp.~21, 23--25]{BickelEtAl2008}, but both lists use modulus $32$ twice. Expanding a repeated class by their ``2 trick'' introduces modulus $64$ \cite[Corollary~4]{BickelEtAl2008}. Since $64$ divides neither listed value, neither is the LCM of the expanded distinct system.

Integer programming has been used for covering-number tests \cite[Appendix~B]{McNewSetty2026} and maximum coverage \cite[Section~6]{HarringtonEtAl2026}. Our construction search introduces increasing blocks and joint revisits, using iterated local search \cite{LourencoEtAl2001} and ruin--recreate perturbations \cite{SchrimpfEtAl2000}. The residue lists establish the bounds independently of these procedures.

\section{Preliminaries and Verification}\label{sec:verification}

\begin{definition}
A \emph{covering system} is a finite collection $\{a_i\pmod{m_i}:1\le i\le n\}$, with $m_i\ge2$ and $0\le a_i<m_i$, such that
\[
 \bigcup_{i=1}^n\{x\in\Z:x\equiv a_i\pmod{m_i}\}=\Z.
\]
The system is \emph{distinct} if its moduli are pairwise distinct. After ordering them as $m_1<\cdots<m_n$, we call $m_1$ its \emph{minimum modulus} and write $L=\lcm(m_1,\ldots,m_n)$.
\end{definition}

Since every modulus divides $L$, the union of the classes is $L$-periodic. Hence coverage of $\Z$ is equivalent to coverage of $0,\ldots,L-1$. The verifier checks the moduli, residues, minimum and actual LCM, then marks $a,a+d,a+2d,\ldots<L$ for each pair $(a,d)$. No point remains unmarked in any of the four lists. Appendix~\ref{app:verifier} gives the verification command; the program uses only integer arithmetic.

\section{The Search Formulation}\label{sec:model}

\subsection{The modulus set}\label{sec:prescribed-moduli}

Let $D=\{d_1,\ldots,d_k\}$ be any nonempty finite set of distinct integers at least $2$, and put
\[
  m=\min D,\qquad L=\lcm(d_1,\ldots,d_k).
\]
The set $D$ is fixed throughout the search; only the residues are changed. Every class has period dividing $L$, so it is enough to check the integers $0,\ldots,L-1$.
For every $d\in D$ choose one residue $a_d\in\Z/d\Z$, and write
\[
  C(d,a_d):=\{x\in\Z/L\Z:x\equiv a_d\pmod d\}.
\]
The maximum-coverage objective, also used in \cite[Section~6]{HarringtonEtAl2026}, is
\begin{equation}\label{eq:objective}
  F_D(a):=\left|\bigcup_{d\in D}C(d,a_d)\right|.
\end{equation}
If $F_D(a)=L$, the selected classes cover all of $\Z$. Otherwise, the output is a partial covering, together with its uncovered set
\[
  H_D(a)=\{x\in\{0,\ldots,L-1\}:x\not\equiv a_d\pmod d
                  \text{ for every }d\in D\}.
\]
Thus $|H_D(a)|=L-F_D(a)$ counts individual uncovered residues, also called \emph{gaps}, within one period.

For the least-common-multiple bounds studied here, we instead begin with a minimum modulus $m$ and a candidate $L$ divisible by $m$, then supply
\[
  D=\D_m(L):=\{d:d\mid L,\ d\ge m\}
\]
to the same search. Using every allowed divisor does not restrict existence: a covering using a subset can be extended by assigning arbitrary residues to the unused divisors. Since $m$ and $L$ are included, the resulting system has minimum modulus $m$ and least common multiple $L$. The experiments below use this full-divisor choice of $D$.

For $m=8$, the separate exclusion computations below reduce the smaller candidate LCMs to two unresolved values.

\subsection[Candidate reduction for m1=8]{Candidate reduction for $m_1=8$}\label{sec:m8-candidates}

We apply the reciprocal-sum and partial-sum filters from our previous joint work \cite[Sections~4--6]{ZhangZhang2026}. Since $8\mid L$ and a construction at $L=20160$ is available, it suffices to consider
\[
 \mathcal R_8=\{L\in\mathbb N:8\mid L,\ 8\le L\le20160\}.
\]
Exact enumeration of the necessary condition $\sum_{d\in\D_8(L)}1/d\ge1$ leaves $53$ candidates. The preliminary exclusions reduce this set to
\[
\begin{aligned}
 \mathcal C_8=\{&10080,12600,15120,15840,16800,17640,\\
                &18480,18720,19800,20160\}.
\end{aligned}
\]
The preliminary exclusions use the earlier lower bound \cite{ZhangZhang2026} and partial-sum tests. The supplement gives the elementary lower-range reductions, including $7920\to5040$, and the computation records. Table~\ref{tab:m8-exclusions} records the subsequent reductions. Depth $k$ means that the filter uses the first $k$ increasing divisors in $\D_8(L)$.

\begin{table}[htbp]
\centering
\caption{Successive exclusions from $\mathcal C_8$ for $m_1=8$. The last column counts all remaining candidates with $L\le20160$, including the constructed value $20160$.}
\label{tab:m8-exclusions}
\begin{tabular}{@{}clc@{}}
\toprule
Depth & Excluded LCMs & Candidates remaining\\
\midrule
23 & $16800,17640$ & 8\\
25 & $18480,18720,19800$ & 5\\
30 & $12600,15840$ & 3\\
\bottomrule
\end{tabular}
\end{table}

Thus the computations leave
\[
 \boxed{L_{\min}(8)\in\{10080,15120,20160\}.}
\]
The value $20160$ is realized by our construction. For $10080$ and $15120$, neither a covering nor an exclusion has been established. Consequently, proving the $m=8$ case of Conjecture~\ref{conj:minimal-lcm} reduces to excluding these two values. The exclusion results use Gurobi's reported bounds and infeasibility decisions; they are distinct from the independently verified covering certificates.

\subsection{Coordinate updates}

Maintain the coverage multiplicity
\[
  c(x):=\#\{d\in D:x\equiv a_d\pmod d\}
\]
for every $x\in\Z/L\Z$. Then
\[
  F_D(a)=\#\{x:c(x)>0\},\qquad
  U_D(a)=L-F_D(a)=\#\{x:c(x)=0\}.
\]
To update the residue associated with a single modulus $d$, first remove its old class and denote the resulting count by $c_{-d}(x)$. The number of newly covered points obtained by assigning residue $r$ is
\begin{equation}\label{eq:best-response}
  g_d(r)=\#\{x:c_{-d}(x)=0,\ x\equiv r\pmod d\}.
\end{equation}
We choose a residue maximizing $g_d(r)$, breaking ties uniformly at random. If $a'$ differs from $a$ only at $d$, with new residue $r$, then
\[
 F_D(a')-F_D(a)=g_d(r)-g_d(a_d).
\]
Thus an exact best response cannot reduce coverage. The counts are updated along the removed and inserted arithmetic progressions.

The count array uses $O(L)$ memory. Computing all gains for one modulus takes $O(L)$ time; adding or removing its class takes $O(L/d)$ time. A sweep through $b$ active moduli costs $O(bL)$. With $T$ workers and $K$ stored uncovered sets, memory use is $O((T+K)L)$, excluding temporary candidate pools. Segmentation reduces the active coordinate count but retains the dependence on $L$.

\subsection{Translation and CRT normalization}

Global translation preserves coverage: adding the same integer $t$ to every residue maps a covering to an equivalent covering. More strongly, let $P\subseteq D$ be a set of pairwise-coprime moduli. The Chinese remainder theorem gives a translation $t$ satisfying
\[
  a_p+t\equiv -1\pmod p\qquad(p\in P).
\]
Consequently, one may fix without loss of generality
\begin{equation}\label{eq:crt-prefix}
  a_p=p-1\qquad(p\in P).
\end{equation}
These fixed classes form the \emph{CRT prefix}. For example, at $L=20160$ we fix residues $8,34,63$ modulo $9,35,64$. The implementation selects a pairwise-coprime subset among the first $h$ increasing moduli, maximizing (product, sum, cardinality) lexicographically by dynamic programming on prime-support masks. Appendix~\ref{app:parameters} lists the prefixes. Their product need not equal $L$; residual translation symmetries may remain.

\section{The Segmented Search}\label{sec:block-search}

\subsection{Insertion and freezing}

After removing the CRT prefix, order the free moduli as
\[
  q_1<q_2<\cdots<q_N.
\]
For a target block length $b$, partition them into consecutive blocks
\[
  \B_1,\B_2,\ldots,\B_s,
  \qquad |\B_j|=b\ (j<s),\qquad s=\lceil N/b\rceil,
\]
with the last block containing the remainder. Alternatively, prescribe a block count $s$ and put $\lceil N/s\rceil$ moduli in each of the first $N\bmod s$ blocks and $\lfloor N/s\rfloor$ in the others. At insertion stage $j$, only $\B_j$ is optimized; the prefix and earlier blocks are fixed, and later blocks are absent. At the end of its trajectory budget the best assignment is retained. Scheduled revisits may then change completed blocks before the next insertion.

If $A_{j-1}$ contains the CRT prefix and the classes fixed in the first $j-1$ blocks, the stage-$j$ subproblem is
\begin{equation}\label{eq:block-objective}
  \max_{\{a_d:d\in\B_j\}}
  \left|
  \bigcup_{(d,a_d)\in A_{j-1}}C(d,a_d)
  \ \cup\!
  \bigcup_{d\in\B_j}C(d,a_d)
  \right|.
\end{equation}
Each stage has fewer active coordinates than the full problem, but is still solved heuristically. Early blocks cannot anticipate every later choice. Revisits allow their residues to change in the presence of classes added subsequently.

\subsection{Searching within a block}

A \emph{trajectory} combines randomized greedy initialization, bounded coordinate search and ruin--recreate perturbations. Initialization samples residues whose gain is within $\max\{1,\lfloor L/(25d)\rfloor\}$ of the maximum. Coordinate search uses exact best responses in shuffled order, with an additional increasing-order sweep after a no-gain round. Perturbations remove and rebuild selected active classes; bounded worsening moves are allowed while a separate historical best is retained. Independent trajectories compete within each insertion budget.

Appendix~\ref{app:trajectory-kernel} gives the stopping, acceptance and perturbation rules. Bounded sweeps do not certify coordinate local optimality: a tied update can enable an improvement at an earlier coordinate. Complete coverage is checked independently of these stopping rules.

\section{Joint Reoptimization}\label{sec:revisit}

\subsection{Joint changes of residues}

The usefulness of a revisit depends on the positions of the uncovered points, as well as their number. Suppose two gaps remain and the only available modulus is $d$. One class modulo $d$ can cover both gaps precisely when they are congruent modulo $d$. Changing an earlier class may move the gaps into a compatible position, even if their number stays the same. A revisit lets the search make such coordinated changes.

A revisit selects a window $W\subseteq\{1,\ldots,s\}$, unfreezes the blocks in $W$, leaves all other blocks fixed, and restarts perturbed coordinate search on
\[
  \bigcup_{j\in W}\B_j.
\]
We use three scales.

\begin{description}
  \item[Adjacent-block revisits.] Reoptimize $\B_j\cup\B_{j+1}$ to repair artificial block boundaries and permit cross-boundary exchanges.
  \item[Subset revisits.] Reoptimize pairs, triples, or larger selections of nonadjacent blocks. These windows address long-range interactions while controlling the number of free coordinates.
  \item[Global revisits.] Unfreeze every free modulus while retaining only the lossless CRT prefix. This activates the largest coordinate set. The historical runs initialize it from a saved assignment; Section~\ref{sec:controlled} compares cold-start procedures.
\end{description}

\subsection{Perturbation schedules and multiple seed assignments}

Seeded revisits start from a saved assignment. Appendix~\ref{app:revisit-schedules} gives their trajectory rules and ruin-size schedules.

For the $m=9$ completion phase, a \emph{seed bank} retains multiple assignments. An \emph{epoch} generates trajectories from the current bank, pools their best states and endpoints with the old bank, and selects the next bank. We use ``seed assignment'' for a saved assignment and ``random seed'' for an integer initializing a pseudorandom generator. 

The bank selects assignments by coverage and by distances between their residue vectors and uncovered sets. Its trajectories allow bounded worsening moves and retain a historical best. Appendix~\ref{app:bank-kernel} gives the exact selection thresholds, acceptance probabilities and schedules. The follow-up in Section~\ref{sec:followup} changes parent sampling and schedule cycling while holding these trajectory rules fixed within its capacity comparison.

\section{Computational Experiments}\label{sec:controlled}

\subsection{The initial experiment}

The initial study uses six configurations and ten random seeds on each of $(m,L)=(8,20160),(8,30240),(9,75600)$, giving $180$ runs. Instances, random seeds and scan limits were fixed before the batch. Every run starts from the CRT prefix, selected from the first $25$ moduli, without an imported assignment. Free-block length $25$ gives partitions $25+25+24$, $25+25+25+11$, and $25+25+25+25+9$.

Each run has a scan limit $Q=240000$, an upper bound on the number of best-response gain scans. One scan evaluates the gain vector over all $L$ points for one modulus; greedy construction and ruin--recreate updates both consume scans. The limit matches full-period scans within an instance, not total operations or elapsed time. The configurations allocate this allowance as follows:
\begin{description}
\item[Whole.] Use the entire allowance for repeated construction trajectories with all free moduli active.
\item[Segmented.] Split the entire allowance equally across insertion blocks, with no revisit.
\item[Adjacent.] Use half for insertions and half, divided equally, for reverse and forward passes over adjacent block pairs. Their revisit schedules alternate large-to-small and mixed perturbations.
\item[Global.] Use half for insertions and half for three full-set revisit phases, with medium, large-to-small and mixed perturbations in that order.
\item[Bank-1 and Bank-20.] Use half for insertions and half, divided equally, for five bank epochs. Both use the same bank trajectory kernel and selection thresholds, with capacity $K=1$ or $20$. Here $K$ is the maximum number of saved seed assignments retained in the bank, not the number actually sampled.
\end{description}
The four half-insertion configurations share insertion random seeds and scan limits. Ordinary revisits restart from a saved assignment; bank trajectories use selected parents. Only the Bank-1/Bank-20 pair isolates capacity; Global/Bank-20 also changes the acceptance kernel.

Runs use one worker, $32$ perturbation steps and at most $80$ coordinate rounds, stopping at the first covering or the scan limit. All final assignments passed independent checks. Execution details and timing are in Appendix~\ref{app:experiment-details}.

\subsection{Initial results}

For each row, let $N$ be the number of runs and $a^{(r)}$ the complete assignment returned by run $r$. We report
\[
 U_r=|H_{D_m(L)}(a^{(r)})|=L-F_{D_m(L)}(a^{(r)}),\qquad
 s=\#\{r:U_r=0\}.
\]
The entry $s/N$ counts successful runs out of all $N$ runs. ``Uncovered residues'' gives the median of $U_1,\ldots,U_N$, followed by their minimum and maximum in brackets, including successful runs with $U_r=0$. The brackets are ranges, not confidence intervals; for even $N$, the median averages the two central values.

For a successful run, let $S_r$ be its cumulative gain-scan count when a covering is first found, including initialization and all preceding stages. ``Scans to success'' reports $\operatorname{median}_{U_r=0}S_r/10^3$, rounded to one decimal place; a dash means $s=0$. This conditional median excludes failed runs. Actual scans may be fewer than the per-run cap $Q$ because of early success. These conventions apply to all experimental tables; CPU timings are reported separately in Table~\ref{tab:controlled-timing}.

\begin{table}[!htb]
\centering
\caption{Initial study: $N=10$ runs per row, each with scan limit $Q=240000$. Uncovered residues are counted in one period $0,\ldots,L-1$ and summarized over all runs; scans to success are summarized over successful runs only. A dash denotes no success.}
\label{tab:controlled-results}
\small
\begin{tabular}{@{}lcrr@{}}
\toprule
Method & \shortstack{Successes\\$s/N$} & \shortstack{Uncovered residues\\median [min, max]} & \shortstack{Scans to success\\median ($10^3$)}\\
\midrule
\multicolumn{4}{c}{$m=8,\ L=20160$}\\
Whole & 0/10 & $221.5\,[122,362]$ & ---\\
Segmented & 0/10 & $60\,[12,102]$ & ---\\
Adjacent & 0/10 & $32.5\,[12,82]$ & ---\\
Global & 3/10 & $33\,[0,79]$ & 129.1\\
Bank-1 & 2/10 & $21\,[0,87]$ & 153.2\\
Bank-20 & 2/10 & $29.5\,[0,91]$ & 153.2\\
\midrule
\multicolumn{4}{c}{$m=8,\ L=30240$}\\
Whole & 0/10 & $127.5\,[85,302]$ & ---\\
Segmented & 8/10 & $0\,[0,50]$ & 180.0\\
Adjacent & 7/10 & $0\,[0,96]$ & 90.0\\
Global & 7/10 & $0\,[0,57]$ & 90.0\\
Bank-1 & 7/10 & $0\,[0,93]$ & 90.0\\
Bank-20 & 8/10 & $0\,[0,93]$ & 90.0\\
\midrule
\multicolumn{4}{c}{$m=9,\ L=75600$}\\
Whole & 0/10 & $1096.5\,[633,1881]$ & ---\\
Segmented & 7/10 & $0\,[0,345]$ & 192.0\\
Adjacent & 6/10 & $0\,[0,423]$ & 96.0\\
Global & 6/10 & $0\,[0,382]$ & 96.0\\
Bank-1 & 6/10 & $0\,[0,397]$ & 96.0\\
Bank-20 & 6/10 & $0\,[0,397]$ & 96.0\\
\bottomrule
\end{tabular}
\end{table}

Whole finds no covering. Segmented succeeds on two instances, while Global finds three coverings at $L=20160$. At $L=30240$ and $75600$, five successes in each half-insertion configuration occur before revisits begin. Bank-20 has no consistent advantage over Bank-1. These results compare complete procedures: segmentation changes initialization order, active coordinates and work allocation together. The follow-up separates initialization order and tests revised bank sampling.

\subsection{The follow-up experiment}\label{sec:followup}

The original bank restarts parent and schedule indices each epoch. Across its Bank-20 runs, only parent positions $1,2$ and two of the three schedules receive work (Appendix~\ref{app:initial-diagnostics}). We therefore designed a separate follow-up after inspecting the initial results, fixing its matrix and random seeds before execution:
\begin{description}
\item[Initialization order: $30$ runs.] On the three original instances, Ordered-Whole changes only greedy initialization from shuffled to increasing order, using the same ten random seeds and $Q=240000$.
\item[Bank sampling: $40$ runs.] On $(8,20160)$, compare $K=1,20$ at $Q=240000,960000$, ten runs per row. Parent positions are shuffled and visited cyclically; schedule counters continue across epochs. Both capacities retain half-insertion and five-epoch allocation.
\item[Larger periods: $60$ runs.] On $(10,302400)$ and $(11,3326400)$, compare Ordered-Whole, Segmented and Global at $Q=15000,60000$, five runs per row. Since each scan costs $O(L)$, work is compared within each period.
\end{description}
The machine, build flags and trajectory kernel are unchanged. All $130$ final assignments passed independent coverage and accounting checks.

\subsection{Follow-up results}

Increasing the initial modulus order reduces the median number of uncovered residues on all three instances, but neither Whole variant finds a covering (Table~\ref{tab:ordering}). Thus ordering alone does not reproduce the Segmented successes.

Table~\ref{tab:revised-bank} shows no consistent advantage for capacity $20$. All forty runs enter bank search; the larger bank reaches parent positions up to $9$ at the lower scan limit and $20$ at the higher one, and uses all three schedules. This removes the original sampling restriction without producing a consistent gain in this sample.

\begin{table}[!htb]
\centering
\caption{Revised bank sampling on $m=8,L=20160$, $N=10$ per row. The scan limit $Q$ is per run; $K$ is the maximum number of retained seed assignments. Uncovered-residue summaries include all runs; scan medians include only successful runs.}
\label{tab:revised-bank}
\small
\begin{tabular}{@{}rccrr@{}}
\toprule
\shortstack{Scan limit\\$Q/10^3$} & \shortstack{Seed capacity\\$K$} & \shortstack{Successes\\$s/N$} & \shortstack{Uncovered residues\\median [min, max]} & \shortstack{Scans to success\\median ($10^3$)}\\
\midrule
240 & 1 & 3/10 & $27\,[0,82]$ & 166.0\\
240 & 20 & 2/10 & $25.5\,[0,78]$ & 149.3\\
960 & 1 & 4/10 & $2.5\,[0,63]$ & 522.8\\
960 & 20 & 4/10 & $4\,[0,73]$ & 522.8\\
\bottomrule
\end{tabular}
\end{table}

In the larger-period tests (Table~\ref{tab:larger-periods}), no method succeeds at $Q=15000$, or at $Q=60000$ for $m=10$. At $m=11,Q=60000$, Ordered-Whole, Segmented and Global succeed in $1/5,5/5,2/5$ runs. Both Global successes occur during insertion, so revisits rescue none of its remaining runs. These small samples do not establish a general ranking. Different total scan limits change stage allocations, so their trajectories are not nested truncations of one run.


\section{Conclusion}\label{sec:scope}

The four residue lists establish the upper bounds in Theorem~\ref{prop:upper-bounds}. For $m=8$, the exclusion computations leave only $L=10080$ and $15120$ below the constructed value $20160$; both remain unresolved. Matching lower bounds in Conjecture~\ref{conj:minimal-lcm} are still open. The experiments support using segmented initialization on the tested instances, while the benefits of revisits and larger seed banks depend on the configuration.

\clearpage
\appendix

\section{Segmented Search Pseudocode}\label{app:pseudocode}

The abstract framework accepts a prescribed set $D$; the recorded implementation constructs $D=\D_m(L)$. Choose the prefix before partitioning, so that $D=P\,\dot\cup\,B_1\,\dot\cup\cdots\dot\cup\,B_s$. An optional complete initial assignment, as used for $m=9$, bypasses fresh insertions and enters the revisit phase. Search budgets bound every optimization call.

\begin{lstlisting}[style=pseudo]
Input: D, prefix-selection and partition rules, budgets,
       revisit schedules, optional complete initial assignment
L <- lcm(D)
P <- pairwise-coprime subset selected from D
fix a[p] <- p-1 for every p in P
partition D minus P into B[1..s]
initialize counts from the prefix

if a complete initial assignment is supplied:
    check its moduli, residues, and prefix; recompute counts
    elite <- initial assignment
else:
    for j <- 1,...,s:
        add B[j]; optimize it within the insertion budget
        keep the best assignment for the completed blocks
        run scheduled revisits of completed blocks
        independently recompute and save the checkpoint
        if coverage is complete:
            assign arbitrary residues to all unused moduli
            break
    elite <- the resulting complete assignment on D

if elite is not a covering and revisit budget remains:
    run configured joint revisits, retaining the historical best
if seed-bank search is enabled and elite is not a covering:
    initialize the bank from elite and any supplied complete seeds
    for each budgeted epoch:
        generate perturbed trajectories; retain best and endpoints
        select the bank by quality, deduplication, and distances
        update elite and save a recomputed checkpoint
        stop if elite is a covering
recompute elite's coverage and uncovered set over 0,...,L-1
return elite and its uncovered set
\end{lstlisting}

\subsection{Insertion trajectory kernel}\label{app:trajectory-kernel}

A \emph{trajectory} is one initialization followed by bounded coordinate search and a specified number of perturbation steps. The retained implementation uses the following rules.

\begin{enumerate}[label=\textbf{S\arabic*.}]
\item \textbf{Randomized greedy initialization.} Visit the active moduli in random order. For modulus $d$, sample uniformly among residues with gain at least $\max_r g_d(r)-\delta_d$, where
\[
 \delta_d=\max\{1,\lfloor L/(25d)\rfloor\}.
\]
Previously unassigned moduli have no old class to remove.
\item \textbf{Coordinate search.} Use exact best responses in a random coordinate order, for at most $80$ rounds. After a round with no gain, try one additional sweep in increasing coordinate order; stop if that sweep also gives no gain. Tie choices remain random in both sweeps. This is a bounded heuristic stopping rule, not a certificate of coordinate local optimality: a later tied update can enable improvement at an earlier coordinate. Every exact best-response update is nevertheless nondecreasing.
\item \textbf{Ruin--recreate perturbations.} For a nonempty active set, select uniformly between $\min(2,b)$ and $\min(6,b)$ coordinates, remove them, rebuild in random order using the same greedy tolerance, and repeat coordinate search. Retain a separate historical best. Accept a new current state if its score does not decrease; on every eighth step, also accept a loss of at most two points relative to the current score. On a rejected twelfth step, restore the historical best. Step counts refer to one-based perturbation indices.
\end{enumerate}

Independent random initializations are attempted within each insertion budget, and the highest-scoring assignment is retained. A sweep or trajectory limit may stop before all coordinate improvements have been exhausted. A local-optimality certificate, if desired, requires checking every coordinate against the same unchanged final state; it is not needed to verify a complete covering.

\subsection{Ordinary revisit schedules}\label{app:revisit-schedules}

Seeded revisits begin with coordinate search on the saved assignment, followed by perturbations. The ordinary revisit kernel uses the acceptance and historical-best rules of Appendix~\ref{app:trajectory-kernel}. Its \emph{large-to-small} schedule linearly decreases a ruin size between specified limits, rounds to an integer, adds a uniform jitter in $\{-1,0,1\}$, and clips to those limits. Its \emph{mixed} schedule chooses a small or large ruin range with probability $1/2$ each and samples uniformly within the selected range. The limits and number of steps are configuration parameters, rather than consequences of the block partition.

\subsection{The smallest LCM case}
For completeness, the classes $0\pmod2$, $1\pmod3$, $3\pmod4$, $5\pmod6$, and $9\pmod{12}$ establish $L_{\min}(2)\le12$. Any smaller candidate is an even integer. For $L=2,4,8,10$, the sum of $1/d$ over divisors $d\ge2$ of $L$ is less than $1$. For $L=6$ the sum equals $1$, but classes modulo $2$ and $3$ necessarily intersect once in a period. Even with a class modulo $6$, their union has size at most $3+2-1+1=5$. No smaller candidate therefore admits a distinct covering.

\section{Exploratory Block-Length Runs}\label{sec:block-size}

Smaller blocks use fewer coordinates in each initial stage, but create more boundaries and possible revisit windows. With $s$ blocks, considering every subset of sizes $2$ through $r$ would involve
\[
 \sum_{i=2}^{\min(r,s)}\binom{s}{i}
\]
windows per enumeration. The recorded runs use specified adjacent or subset schedules. Selecting windows from measured compatibility with uncovered points is a possible extension, not a procedure evaluated here.

\subsection{Recorded runs for \texorpdfstring{$m=8$}{m=8}}

Table~\ref{tab:15120-block-count} lists single saved outcomes for $m=8,L=15120$, with prefix $16,27,35$ and $70$ free moduli. The recorded insertion allowance is $1200$ trajectories per worker and $18$ perturbation steps, with ten workers; each scheduled revisit has $160$ trajectories per worker and $18$ steps. The three- to five-block runs include all pairs and triples, whereas the six-block run includes subset sizes two through six. The one- and two-block launch script enables subset refinement: only the two-block case has a subset of size at least two. These schedules entail different amounts of work despite matching per-call allowances. All six checkpoints assign all $73$ allowed moduli; none is an exclusion certificate.

\begin{table}[htbp]
\centering
\caption{Individual saved outcomes at $m=8,L=15120$. These are not equal-budget estimates of block-count effects.}
\label{tab:15120-block-count}
\begin{tabular}{@{}clrr@{}}
\toprule
Blocks & Actual block sizes & Covered & Uncovered\\
\midrule
1 & $70$ & 15019 & 101\\
2 & $35+35$ & 15026 & 94\\
3 & $24+23+23$ & 15043 & 77\\
4 & $18+18+17+17$ & 15027 & 93\\
5 & $14+14+14+14+14$ & 15052 & 68\\
6 & $12+12+12+12+11+11$ & 15027 & 93\\
\bottomrule
\end{tabular}
\end{table}

The five-block run has the fewest uncovered points among these saved assignments. This observation alone does not identify an optimal block length or a generally superior schedule.

Table~\ref{tab:target-size} gives individual outcomes for fixed target block lengths at two other periods. Every listed $L=30240$ run completed a covering; among the listed $L=27720$ runs, target length $25$ left the fewest points uncovered. These observations do not estimate success probabilities or sensitivity across independent repetitions.

\begin{table}[htbp]
\centering
\caption{Saved outcomes for fixed target lengths, with $m=8$. Each final assignment contains all $89$ allowed moduli.}
\label{tab:target-size}
\small
\begin{tabular}{@{}crrrr@{}}
\toprule
Target length & \multicolumn{2}{c}{$L=27720$} & \multicolumn{2}{c}{$L=30240$}\\
\cmidrule(lr){2-3}\cmidrule(l){4-5}
 & Covered & Uncovered & Covered & Uncovered\\
\midrule
15 & 27148 & 572 & 30240 & 0\\
17 & 27133 & 587 & 30240 & 0\\
20 & 27147 & 573 & 30240 & 0\\
22 & 27147 & 573 & 30240 & 0\\
25 & 27205 & 515 & 30240 & 0\\
\bottomrule
\end{tabular}
\end{table}

At $L=20160$, a saved target-length-$15$ state covers $20083$ points but contains only the prefix and four blocks of $15$ free moduli, or $63$ classes in total. The last $14$ free moduli have not been inserted. It is therefore an incomplete construction state and is not compared with the complete assignment as a block-length benchmark.

\section{Recorded Construction Paths}\label{sec:records}

The following paths document successful historical searches. They are not controlled comparisons and are not needed for the certificate proofs. ``After insertion'' denotes the state before its subsequent revisits; ``after revisits'' denotes the state used for the next insertion. Final certificate restrictions need not reproduce earlier states, since revisits can change their residues. Appendix~\ref{app:parameters} summarizes the recovered configurations.

\subsection{\texorpdfstring{$m=8,L=20160$}{m=8,L=20160}: three blocks and a cross-block revisit}

The CRT prefix is
\[
  8\pmod9,\qquad 34\pmod{35},\qquad 63\pmod{64},
\]
where $9,35,64$ are pairwise coprime and have product $20160$. The remaining $74$ moduli are divided into blocks of sizes $25+25+24$. Table~\ref{tab:m8-stages} distinguishes the coverage immediately after an insertion from the coverage after a revisit.

\begin{table}[htbp]
\centering
\caption{Key stages for $m=8,L=20160$.}
\label{tab:m8-stages}
\small\begin{tabularx}{\textwidth}{@{}lrrX@{}}
\toprule
Stage & Covered & Uncovered & Description\\
\midrule
CRT prefix $+$ block 1 & 18439 & 1721 & first block contains 25 smaller free moduli\\
Insert block 2 & 19923 & 237 & block 1 frozen while block 2 is optimized\\
Revisit blocks 1 and 2 & 19962 & 198 & jointly reopen 50 moduli; gain 39 points\\
Insert block 3 & 20160 & 0 & final 24 moduli complete the cover\\
\bottomrule
\end{tabularx}
\end{table}

In this recorded path, reopening the first two blocks increases coverage by $39$ points, and the following insertion completes the covering. This trace does not show that a comparable search without that revisit would fail.

\subsection{\texorpdfstring{$m=9,L=75600$}{m=9,L=75600}: projected seeds and joint revisits}

The CRT prefix is $25,27,56$, and the $109$ free moduli are partitioned as $19+18+18+18+18+18$. This run resumes a full assignment projected from a saved $L=151200$ construction. The projection retains classes whose moduli divide $75600$, translates the inherited residues to agree with the target prefix, and greedily fills missing moduli. It then performs up to three sweeps on newly filled coordinates. A \emph{preserve} variant leaves other inherited coordinates fixed; a \emph{relaxed} variant additionally performs up to three sweeps on all free coordinates. A sweep with no gain ends either sweep loop early.

The saved projection summary records $87$ inherited and $25$ missing moduli and translation $0$ modulo $37800$. The preserve and relaxed checkpoints cover $74104$ and $75077$ points, respectively. The successful run starts from the relaxed checkpoint and also supplies the preserve checkpoint to the later seed bank. Thus its initialization includes information from an earlier construction. The supplementary \texttt{projection\_m9} directory now includes the source $151200$ assignment, projection program, both target checkpoints, and a reconstructed command that reproduces both checkpoints byte for byte. This establishes the saved-data chain, not the original executable build or total upstream search cost. In contrast, the initial experiment in Section~\ref{sec:controlled} finds covers in seven of ten fresh Segmented runs on this instance without projection or a bank.

Two rounds of adjacent revisits bring coverage to $75173$. Global revisits then give
\[
 75173\longrightarrow75558\longrightarrow75594.
\]
The bank starts with two seeds. Its target capacity is $20$, reached at the end of the first epoch. The best numbers of uncovered points after epochs $1,2,3$ are $6,4,0$, respectively. Each epoch has a nominal allowance of $10000$ trajectories, but the third stops early when a covering is found. The historical file records $30000$ by multiplying the epoch index by its allowance; this is not an actual executed-trajectory count. The exact total is not recovered, and no comparison of runtime or success probability is drawn from this run.

\subsection{\texorpdfstring{$m=10,L=302400$}{m=10,L=302400}: nine balanced blocks}

The prefix is $25,32,63$; the $156$ free moduli form three blocks of $18$ and six of $17$. Table~\ref{tab:insertion-revisit} separates insertions from the intervening adjacent revisits. In particular, insertion of block $8$ leaves $101$ points uncovered, its revisit reduces that number to $84$, and insertion of block $9$ completes the covering.

\subsection{\texorpdfstring{$m=11,L=3326400$}{m=11,L=3326400}: completion after block 11}

The prefix is $11,27,32,35$. The $322$ free moduli are partitioned into sixteen blocks of $18$ and two of $17$, so eighteen blocks are planned. After insertion of block $11$, four points remain uncovered. Reoptimizing blocks $10$ and $11$ then covers all four. This is the first complete covering in the saved path; block $12$ does not cause completion. Blocks $12$--$18$ are subsequently added to retain the full-divisor output with $326$ classes.

\begin{table}[H]
\centering
\caption{Coverage immediately after each insertion and after its subsequent adjacent revisits. A dash denotes a stage beyond the end of the $m=10$ construction. The $m=11$ covering is already complete after the block-$11$ revisit.}
\label{tab:insertion-revisit}
\begin{tabular}{@{}crrrr@{}}
\toprule
 & \multicolumn{2}{c}{$m=10$, $L=302400$} & \multicolumn{2}{c}{$m=11$, $L=3326400$}\\
\cmidrule(lr){2-3}\cmidrule(l){4-5}
Block & Insertion & Revisits & Insertion & Revisits\\
\midrule
1 & 228224 & 228224 & 2473986 & 2473986\\
2 & 270257 & 273197 & 2922924 & 2922999\\
3 & 288733 & 289046 & 3114358 & 3117758\\
4 & 295788 & 295921 & 3209512 & 3210660\\
5 & 298942 & 299151 & 3256513 & 3259369\\
6 & 300968 & 300998 & 3293305 & 3294892\\
7 & 301886 & 301894 & 3308504 & 3309727\\
8 & 302299 & 302316 & 3317812 & 3318872\\
9 & 302400 & 302400 & 3322637 & 3322922\\
10 & --- & --- & 3325751 & 3325788\\
11 & --- & --- & 3326396 & 3326400\\
\bottomrule
\end{tabular}
\end{table}

\section{Recorded Configurations}\label{app:parameters}

Table~\ref{tab:parameters} gives partitions and completion operations. A block partition alone is not a complete run specification. The detailed rules in the text describe the retained implementation; a source snapshot and historical executable have not been tied to a verified original build. The provenance file distinguishes recorded settings from missing information.

\begin{table}[H]
\centering
\caption{Partitions and completion operations of the four recorded constructions.}
\label{tab:parameters}
\small
\begin{tabularx}{\textwidth}{@{}c@{\hspace{1.1em}}r@{\hspace{1.1em}}llX@{}}
\toprule
$m$ & $L$ & Prefix moduli & Free-block sizes & Completion operation\\
\midrule
8 & 20160 & $9,35,64$ & $25+25+24$ & first-two-block revisit, then third insertion\\
9 & 75600 & $25,27,56$ & $19+5\times18$ & projected seeds, global revisits, third bank epoch\\
10 & 302400 & $25,32,63$ & $3\times18+6\times17$ & ninth insertion, following adjacent revisits\\
11 & 3326400 & $11,27,32,35$ & $16\times18+2\times17$ & revisit of blocks $10$ and $11$\\
\bottomrule
\end{tabularx}
\end{table}

For $m=9$, the retained launch script specifies prefix depth $25$, ten workers, $200$ trajectories per worker for each revisit schedule, $32$ perturbation steps, two reverse/forward adjacent-revisit rounds, and final global refinement. The bank capacity is $20$, with a nominal $500$ trajectories per target seed slot per epoch and at most four epochs. It stops after three. The base random seed is $20260814075600$. The initial bank has two supplied states, so its first epoch allocates the nominal $10000$ trajectories across those states rather than across twenty pre-existing seeds.

For the bank's $109$ free coordinates, trajectory indices cycle through three schedules. The medium schedule samples a ruin size uniformly from $7$ to $24$. The large-to-small schedule decreases from $54$ to $2$, with rounding and uniform jitter in $\{-2,-1,0,1,2\}$, clipped to that range. The mixed schedule chooses a uniform size from $27$ to $54$ with probability $0.6$, otherwise from $2$ to $7$. All use $32$ perturbation steps and the bank acceptance rule in Section~\ref{sec:revisit}. These are implementation settings, not fitted success-rate estimates.

The $m=11$ successful run belongs to the extended candidate queue. Its log records $200$ nominal trajectories for each insertion and $50$ per adjacent revisit schedule across all workers. The matching script specifies ten workers, $32$ perturbation steps for insertions and revisits, prefix depth $20$, and base seed $20260815110000+3326400$. It reaches full coverage during the first trajectory reported for the block-$10$--$11$ large-to-small revisit. These local allowances do not supply total discovery cost. Hardware details and actual stopping counts are not filled from defaults of unrelated launch scripts.

\subsection{Seed-bank trajectory and selection rules}\label{app:bank-kernel}

The following rules define the retained bank kernel used in the historical $m=9$ completion and both experimental batches. Parent sampling is a separate rule: the initial batch restarts its parent and schedule indices each epoch, whereas the follow-up shuffles parent positions and carries its schedule index across epochs.

For assignments $a,a'$ on the same free moduli, define
\[
 d_R(a,a')=\#\{d:a_d\ne a'_d\},\qquad
 d_H(a,a')=|H_D(a)\triangle H_D(a')|.
\]
Candidate seed assignments are ordered by decreasing coverage, then lexicographically by their uncovered sets and residue vectors. After removing duplicate residue vectors and candidates more than $15$ points below the best, reject a candidate if some retained assignment satisfies both $d_H<8$ and $d_R<12$. Stop at the configured capacity: $20$ in the historical run, and $1$ or $20$ in the experiments. If fewer survive, fill the remaining positions from the same quality-filtered pool of distinct vectors, irrespective of the distance thresholds.

The bank's trajectory kernel differs from ordinary revisits. For a proposed score $v<v_{\rm cur}$, it permits acceptance only if $v\ge v_{\rm best}-12$, and then with probability
\[
 \exp((v-v_{\rm cur})/\theta_t),\qquad
 \theta_t=\max\{0.75,12(1-0.75t)\},
\]
where $t$ is fractional progress through the perturbation steps. A rejected tenth step restores the historical best. Nondecreasing moves are accepted, and the best score is stored separately. Appendix~\ref{app:parameters} specifies the ruin ranges used for this run.

\section{Complete Constructions}\label{app:constructions}

Each entry $a\,(d)$ denotes the residue class $x\equiv a\pmod d$.
Entries are ordered by increasing modulus, read from left to right across each row.
The four lists below preserve the full-divisor assignments produced by the searches.
They contain all data needed to verify Theorem~\ref{prop:upper-bounds}.

\subsection{Minimum modulus 8}
The least common multiple is $L=20160$.

\begingroup
\footnotesize
\setlength{\tabcolsep}{7pt}
\renewcommand{\arraystretch}{1.2}
\begin{longtable}{@{}rrrrr@{}}
\caption{Complete covering for $m=8,L=20160$.}\label{tab:cover-8}\\
\toprule
\certhead
\midrule
\endfirsthead
\multicolumn{5}{c}{Table \thetable\ (continued): $m=8,L=20160$}\\
\toprule
\certhead
\midrule
\endhead
\midrule
\multicolumn{5}{r}{Continued on the next page}\\
\endfoot
\bottomrule
\endlastfoot
\coverentry{3}{8} & \coverentry{8}{9} & \coverentry{6}{10} & \coverentry{1}{12} & \coverentry{2}{14} \\
\coverentry{3}{15} & \coverentry{7}{16} & \coverentry{2}{18} & \coverentry{0}{20} & \coverentry{19}{21} \\
\coverentry{9}{24} & \coverentry{10}{28} & \coverentry{12}{30} & \coverentry{15}{32} & \coverentry{34}{35} \\
\coverentry{29}{36} & \coverentry{10}{40} & \coverentry{22}{42} & \coverentry{23}{45} & \coverentry{21}{48} \\
\coverentry{46}{56} & \coverentry{24}{60} & \coverentry{28}{63} & \coverentry{63}{64} & \coverentry{14}{70} \\
\coverentry{41}{72} & \coverentry{77}{80} & \coverentry{4}{84} & \coverentry{32}{90} & \coverentry{45}{96} \\
\coverentry{13}{105} & \coverentry{52}{112} & \coverentry{30}{120} & \coverentry{112}{126} & \coverentry{54}{140} \\
\coverentry{5}{144} & \coverentry{61}{160} & \coverentry{136}{168} & \coverentry{104}{180} & \coverentry{159}{192} \\
\coverentry{202}{210} & \coverentry{95}{224} & \coverentry{189}{240} & \coverentry{70}{252} & \coverentry{74}{280} \\
\coverentry{77}{288} & \coverentry{239}{315} & \coverentry{287}{320} & \coverentry{220}{336} & \coverentry{230}{360} \\
\coverentry{414}{420} & \coverentry{415}{448} & \coverentry{285}{480} & \coverentry{466}{504} & \coverentry{18}{560} \\
\coverentry{479}{576} & \coverentry{448}{630} & \coverentry{127}{672} & \coverentry{365}{720} & \coverentry{790}{840} \\
\coverentry{95}{960} & \coverentry{298}{1008} & \coverentry{895}{1120} & \coverentry{952}{1260} & \coverentry{31}{1344} \\
\coverentry{509}{1440} & \coverentry{298}{1680} & \coverentry{1183}{2016} & \coverentry{671}{2240} & \coverentry{1642}{2520} \\
\coverentry{2399}{2880} & \coverentry{1631}{3360} & \coverentry{2527}{4032} & \coverentry{802}{5040} & \coverentry{2591}{6720} \\
\coverentry{7871}{10080} & \coverentry{6431}{20160} &  &  &  \\
\end{longtable}
\endgroup

\subsection{Minimum modulus 9}
The least common multiple is $L=75600$.

\begingroup
\footnotesize
\setlength{\tabcolsep}{7pt}
\renewcommand{\arraystretch}{1.2}
\begin{longtable}{@{}rrrrr@{}}
\caption{Complete covering for $m=9,L=75600$.}\label{tab:cover-9}\\
\toprule
\certhead
\midrule
\endfirsthead
\multicolumn{5}{c}{Table \thetable\ (continued): $m=9,L=75600$}\\
\toprule
\certhead
\midrule
\endhead
\midrule
\multicolumn{5}{r}{Continued on the next page}\\
\endfoot
\bottomrule
\endlastfoot
\coverentry{2}{9} & \coverentry{0}{10} & \coverentry{1}{12} & \coverentry{7}{14} & \coverentry{6}{15} \\
\coverentry{3}{16} & \coverentry{14}{18} & \coverentry{12}{20} & \coverentry{15}{21} & \coverentry{7}{24} \\
\coverentry{24}{25} & \coverentry{26}{27} & \coverentry{13}{28} & \coverentry{16}{30} & \coverentry{3}{35} \\
\coverentry{5}{36} & \coverentry{2}{40} & \coverentry{9}{42} & \coverentry{13}{45} & \coverentry{43}{48} \\
\coverentry{4}{50} & \coverentry{44}{54} & \coverentry{55}{56} & \coverentry{48}{60} & \coverentry{54}{63} \\
\coverentry{45}{70} & \coverentry{23}{72} & \coverentry{19}{75} & \coverentry{62}{80} & \coverentry{81}{84} \\
\coverentry{88}{90} & \coverentry{64}{100} & \coverentry{24}{105} & \coverentry{17}{108} & \coverentry{27}{112} \\
\coverentry{22}{120} & \coverentry{33}{126} & \coverentry{116}{135} & \coverentry{18}{140} & \coverentry{59}{144} \\
\coverentry{84}{150} & \coverentry{39}{168} & \coverentry{69}{175} & \coverentry{28}{180} & \coverentry{8}{189} \\
\coverentry{114}{200} & \coverentry{87}{210} & \coverentry{71}{216} & \coverentry{184}{225} & \coverentry{102}{240} \\
\coverentry{201}{252} & \coverentry{8}{270} & \coverentry{98}{280} & \coverentry{34}{300} & \coverentry{138}{315} \\
\coverentry{123}{336} & \coverentry{294}{350} & \coverentry{118}{360} & \coverentry{89}{378} & \coverentry{214}{400} \\
\coverentry{258}{420} & \coverentry{395}{432} & \coverentry{34}{450} & \coverentry{327}{504} & \coverentry{219}{525} \\
\coverentry{89}{540} & \coverentry{58}{560} & \coverentry{214}{600} & \coverentry{305}{630} & \coverentry{494}{675} \\
\coverentry{58}{700} & \coverentry{658}{720} & \coverentry{737}{756} & \coverentry{798}{840} & \coverentry{784}{900} \\
\coverentry{332}{945} & \coverentry{75}{1008} & \coverentry{144}{1050} & \coverentry{858}{1080} & \coverentry{414}{1200} \\
\coverentry{978}{1260} & \coverentry{1034}{1350} & \coverentry{198}{1400} & \coverentry{359}{1512} & \coverentry{1419}{1575} \\
\coverentry{538}{1680} & \coverentry{1278}{1800} & \coverentry{683}{1890} & \coverentry{618}{2100} & \coverentry{1142}{2160} \\
\coverentry{2458}{2520} & \coverentry{359}{2700} & \coverentry{338}{2800} & \coverentry{1115}{3024} & \coverentry{894}{3150} \\
\coverentry{978}{3600} & \coverentry{1817}{3780} & \coverentry{1038}{4200} & \coverentry{4679}{4725} & \coverentry{4618}{5040} \\
\coverentry{4814}{5400} & \coverentry{6078}{6300} & \coverentry{7487}{7560} & \coverentry{6498}{8400} & \coverentry{1439}{9450} \\
\coverentry{1979}{10800} & \coverentry{1878}{12600} & \coverentry{3707}{15120} & \coverentry{12239}{18900} & \coverentry{16578}{25200} \\
\coverentry{31679}{37800} & \coverentry{25423}{75600} &  &  &  \\
\end{longtable}
\endgroup

\newpage
\subsection{Minimum modulus 10}
The least common multiple is $L=302400$.

\begingroup
\footnotesize
\setlength{\tabcolsep}{7pt}
\renewcommand{\arraystretch}{1.15}
\begin{longtable}{@{}rrrrr@{}}
\caption{Complete covering for $m=10,L=302400$.}\label{tab:cover-10}\\
\toprule
\certhead
\midrule
\endfirsthead
\multicolumn{5}{c}{Table \thetable\ (continued): $m=10,L=302400$}\\
\toprule
\certhead
\midrule
\endhead
\midrule
\multicolumn{5}{r}{Continued on the next page}\\
\endfoot
\bottomrule
\endlastfoot
\coverentry{8}{10} & \coverentry{9}{12} & \coverentry{4}{14} & \coverentry{5}{15} & \coverentry{7}{16} \\
\coverentry{13}{18} & \coverentry{2}{20} & \coverentry{17}{21} & \coverentry{3}{24} & \coverentry{24}{25} \\
\coverentry{25}{27} & \coverentry{8}{28} & \coverentry{0}{30} & \coverentry{31}{32} & \coverentry{21}{35} \\
\coverentry{1}{36} & \coverentry{32}{40} & \coverentry{5}{42} & \coverentry{10}{45} & \coverentry{35}{48} \\
\coverentry{14}{50} & \coverentry{43}{54} & \coverentry{20}{56} & \coverentry{6}{60} & \coverentry{62}{63} \\
\coverentry{47}{64} & \coverentry{26}{70} & \coverentry{19}{72} & \coverentry{59}{75} & \coverentry{12}{80} \\
\coverentry{29}{84} & \coverentry{40}{90} & \coverentry{11}{96} & \coverentry{54}{100} & \coverentry{86}{105} \\
\coverentry{61}{108} & \coverentry{100}{112} & \coverentry{96}{120} & \coverentry{41}{126} & \coverentry{115}{135} \\
\coverentry{106}{140} & \coverentry{125}{144} & \coverentry{84}{150} & \coverentry{132}{160} & \coverentry{137}{168} \\
\coverentry{119}{175} & \coverentry{173}{180} & \coverentry{115}{189} & \coverentry{15}{192} & \coverentry{104}{200} \\
\coverentry{16}{210} & \coverentry{115}{216} & \coverentry{123}{224} & \coverentry{109}{225} & \coverentry{17}{240} \\
\coverentry{77}{252} & \coverentry{70}{270} & \coverentry{216}{280} & \coverentry{155}{288} & \coverentry{294}{300} \\
\coverentry{146}{315} & \coverentry{212}{320} & \coverentry{53}{336} & \coverentry{194}{350} & \coverentry{317}{360} \\
\coverentry{178}{378} & \coverentry{204}{400} & \coverentry{286}{420} & \coverentry{341}{432} & \coverentry{180}{448} \\
\coverentry{184}{450} & \coverentry{143}{480} & \coverentry{149}{504} & \coverentry{44}{525} & \coverentry{413}{540} \\
\coverentry{156}{560} & \coverentry{271}{576} & \coverentry{144}{600} & \coverentry{346}{630} & \coverentry{155}{672} \\
\coverentry{259}{675} & \coverentry{694}{700} & \coverentry{137}{720} & \coverentry{10}{756} & \coverentry{4}{800} \\
\coverentry{276}{840} & \coverentry{485}{864} & \coverentry{29}{900} & \coverentry{451}{945} & \coverentry{527}{960} \\
\coverentry{221}{1008} & \coverentry{394}{1050} & \coverentry{593}{1080} & \coverentry{532}{1120} & \coverentry{689}{1200} \\
\coverentry{986}{1260} & \coverentry{1103}{1344} & \coverentry{34}{1350} & \coverentry{944}{1400} & \coverentry{923}{1440} \\
\coverentry{665}{1512} & \coverentry{1469}{1575} & \coverentry{404}{1600} & \coverentry{1577}{1680} & \coverentry{655}{1728} \\
\coverentry{389}{1800} & \coverentry{136}{1890} & \coverentry{1115}{2016} & \coverentry{1294}{2100} & \coverentry{197}{2160} \\
\coverentry{52}{2240} & \coverentry{239}{2400} & \coverentry{1913}{2520} & \coverentry{629}{2700} & \coverentry{44}{2800} \\
\coverentry{1403}{2880} & \coverentry{1673}{3024} & \coverentry{94}{3150} & \coverentry{3227}{3360} & \coverentry{1289}{3600} \\
\coverentry{2656}{3780} & \coverentry{539}{4032} & \coverentry{444}{4200} & \coverentry{3977}{4320} & \coverentry{4294}{4725} \\
\coverentry{1679}{4800} & \coverentry{1901}{5040} & \coverentry{5129}{5400} & \coverentry{5204}{5600} & \coverentry{917}{6048} \\
\coverentry{1994}{6300} & \coverentry{911}{6720} & \coverentry{2939}{7200} & \coverentry{233}{7560} & \coverentry{7289}{8400} \\
\coverentry{6137}{8640} & \coverentry{7444}{9450} & \coverentry{1787}{10080} & \coverentry{2789}{10800} & \coverentry{1204}{11200} \\
\coverentry{9557}{12096} & \coverentry{3929}{12600} & \coverentry{629}{14400} & \coverentry{7638}{15120} & \coverentry{13979}{16800} \\
\coverentry{1144}{18900} & \coverentry{11573}{20160} & \coverentry{16119}{21600} & \coverentry{23069}{25200} & \coverentry{9557}{30240} \\
\coverentry{10319}{33600} & \coverentry{12329}{37800} & \coverentry{22247}{43200} & \coverentry{25979}{50400} & \coverentry{45017}{60480} \\
\coverentry{54089}{75600} & \coverentry{50939}{100800} & \coverentry{47657}{151200} & \coverentry{63989}{302400} &  \\
\end{longtable}
\endgroup

\subsection{Minimum modulus 11}
The least common multiple is $L=3326400$.

\begingroup
\footnotesize
\setlength{\tabcolsep}{7pt}
\renewcommand{\arraystretch}{1.2}
\begin{longtable}{@{}rrrrr@{}}
\caption{Complete covering for $m=11,L=3326400$.}\label{tab:cover-11}\\
\toprule
\certhead
\midrule
\endfirsthead
\multicolumn{5}{c}{Table \thetable\ (continued): $m=11,L=3326400$}\\
\toprule
\certhead
\midrule
\endhead
\midrule
\multicolumn{5}{r}{Continued on the next page}\\
\endfoot
\bottomrule
\endlastfoot
\coverentry{10}{11} & \coverentry{9}{12} & \coverentry{12}{14} & \coverentry{10}{15} & \coverentry{7}{16} \\
\coverentry{11}{18} & \coverentry{16}{20} & \coverentry{4}{21} & \coverentry{6}{22} & \coverentry{3}{24} \\
\coverentry{17}{25} & \coverentry{26}{27} & \coverentry{2}{28} & \coverentry{1}{30} & \coverentry{31}{32} \\
\coverentry{14}{33} & \coverentry{34}{35} & \coverentry{5}{36} & \coverentry{8}{40} & \coverentry{18}{42} \\
\coverentry{2}{44} & \coverentry{19}{45} & \coverentry{36}{48} & \coverentry{2}{50} & \coverentry{17}{54} \\
\coverentry{9}{55} & \coverentry{22}{56} & \coverentry{13}{60} & \coverentry{7}{63} & \coverentry{28}{64} \\
\coverentry{8}{66} & \coverentry{48}{70} & \coverentry{59}{72} & \coverentry{7}{75} & \coverentry{45}{77} \\
\coverentry{64}{80} & \coverentry{42}{84} & \coverentry{18}{88} & \coverentry{49}{90} & \coverentry{15}{96} \\
\coverentry{11}{99} & \coverentry{12}{100} & \coverentry{73}{105} & \coverentry{89}{108} & \coverentry{20}{110} \\
\coverentry{106}{112} & \coverentry{43}{120} & \coverentry{112}{126} & \coverentry{38}{132} & \coverentry{34}{135} \\
\coverentry{6}{140} & \coverentry{95}{144} & \coverentry{127}{150} & \coverentry{136}{154} & \coverentry{108}{160} \\
\coverentry{125}{165} & \coverentry{150}{168} & \coverentry{97}{175} & \coverentry{84}{176} & \coverentry{0}{180} \\
\coverentry{106}{189} & \coverentry{60}{192} & \coverentry{176}{198} & \coverentry{32}{200} & \coverentry{90}{210} \\
\coverentry{35}{216} & \coverentry{60}{220} & \coverentry{162}{224} & \coverentry{172}{225} & \coverentry{158}{231} \\
\coverentry{24}{240} & \coverentry{154}{252} & \coverentry{62}{264} & \coverentry{259}{270} & \coverentry{97}{275} \\
\coverentry{268}{280} & \coverentry{204}{288} & \coverentry{143}{297} & \coverentry{37}{300} & \coverentry{210}{308} \\
\coverentry{304}{315} & \coverentry{204}{320} & \coverentry{320}{330} & \coverentry{234}{336} & \coverentry{62}{350} \\
\coverentry{172}{352} & \coverentry{240}{360} & \coverentry{187}{378} & \coverentry{224}{385} & \coverentry{86}{396} \\
\coverentry{272}{400} & \coverentry{346}{420} & \coverentry{143}{432} & \coverentry{40}{440} & \coverentry{274}{448} \\
\coverentry{247}{450} & \coverentry{452}{462} & \coverentry{463}{480} & \coverentry{218}{495} & \coverentry{436}{504} \\
\coverentry{262}{525} & \coverentry{68}{528} & \coverentry{300}{540} & \coverentry{522}{550} & \coverentry{184}{560} \\
\coverentry{498}{576} & \coverentry{44}{594} & \coverentry{187}{600} & \coverentry{92}{616} & \coverentry{94}{630} \\
\coverentry{200}{660} & \coverentry{66}{672} & \coverentry{547}{675} & \coverentry{56}{693} & \coverentry{482}{700} \\
\coverentry{700}{704} & \coverentry{104}{720} & \coverentry{232}{756} & \coverentry{114}{770} & \coverentry{716}{792} \\
\coverentry{332}{800} & \coverentry{482}{825} & \coverentry{148}{840} & \coverentry{696}{864} & \coverentry{148}{880} \\
\coverentry{72}{900} & \coverentry{686}{924} & \coverentry{79}{945} & \coverentry{684}{960} & \coverentry{416}{990} \\
\coverentry{738}{1008} & \coverentry{472}{1050} & \coverentry{332}{1056} & \coverentry{120}{1080} & \coverentry{262}{1100} \\
\coverentry{946}{1120} & \coverentry{254}{1155} & \coverentry{1142}{1188} & \coverentry{1087}{1200} & \coverentry{386}{1232} \\
\coverentry{904}{1260} & \coverentry{764}{1320} & \coverentry{498}{1344} & \coverentry{97}{1350} & \coverentry{980}{1386} \\
\coverentry{1072}{1400} & \coverentry{300}{1440} & \coverentry{1043}{1485} & \coverentry{66}{1512} & \coverentry{1148}{1540} \\
\coverentry{772}{1575} & \coverentry{200}{1584} & \coverentry{72}{1600} & \coverentry{632}{1650} & \coverentry{344}{1680} \\
\coverentry{1128}{1728} & \coverentry{1468}{1760} & \coverentry{1372}{1800} & \coverentry{1772}{1848} & \coverentry{394}{1890} \\
\coverentry{422}{1925} & \coverentry{1868}{1980} & \coverentry{402}{2016} & \coverentry{3}{2079} & \coverentry{2072}{2100} \\
\coverentry{1916}{2112} & \coverentry{480}{2160} & \coverentry{532}{2200} & \coverentry{498}{2240} & \coverentry{2174}{2310} \\
\coverentry{548}{2376} & \coverentry{562}{2400} & \coverentry{1170}{2464} & \coverentry{1672}{2475} & \coverentry{28}{2520} \\
\coverentry{104}{2640} & \coverentry{2172}{2700} & \coverentry{2350}{2772} & \coverentry{72}{2800} & \coverentry{562}{2880} \\
\coverentry{746}{2970} & \coverentry{24}{3024} & \coverentry{1624}{3080} & \coverentry{997}{3150} & \coverentry{1042}{3168} \\
\coverentry{2372}{3300} & \coverentry{274}{3360} & \coverentry{1767}{3465} & \coverentry{1202}{3520} & \coverentry{2172}{3600} \\
\coverentry{920}{3696} & \coverentry{1032}{3780} & \coverentry{1882}{3850} & \coverentry{3104}{3960} & \coverentry{3090}{4032} \\
\coverentry{3594}{4158} & \coverentry{672}{4200} & \coverentry{3720}{4320} & \coverentry{4072}{4400} & \coverentry{4412}{4620} \\
\coverentry{4497}{4725} & \coverentry{3432}{4752} & \coverentry{2099}{4800} & \coverentry{2034}{4928} & \coverentry{4092}{4950} \\
\coverentry{816}{5040} & \coverentry{3773}{5280} & \coverentry{4438}{5400} & \coverentry{2880}{5544} & \coverentry{4862}{5600} \\
\coverentry{336}{5775} & \coverentry{3146}{5940} & \coverentry{5041}{6048} & \coverentry{1904}{6160} & \coverentry{4036}{6300} \\
\coverentry{2721}{6336} & \coverentry{2556}{6600} & \coverentry{1010}{6720} & \coverentry{5897}{6930} & \coverentry{886}{7200} \\
\coverentry{7261}{7392} & \coverentry{5605}{7425} & \coverentry{3130}{7560} & \coverentry{1970}{7700} & \coverentry{774}{7920} \\
\coverentry{6649}{8316} & \coverentry{5839}{8400} & \coverentry{5981}{8640} & \coverentry{4689}{8800} & \coverentry{4243}{9240} \\
\coverentry{7730}{9450} & \coverentry{1719}{9504} & \coverentry{7105}{9900} & \coverentry{6726}{10080} & \coverentry{6187}{10395} \\
\coverentry{8945}{10560} & \coverentry{4943}{10800} & \coverentry{3786}{11088} & \coverentry{5044}{11200} & \coverentry{3673}{11550} \\
\coverentry{9926}{11880} & \coverentry{10181}{12096} & \coverentry{1274}{12320} & \coverentry{8343}{12600} & \coverentry{10026}{13200} \\
\coverentry{11411}{13860} & \coverentry{2011}{14400} & \coverentry{1552}{14784} & \coverentry{8389}{14850} & \coverentry{2510}{15120} \\
\coverentry{15235}{15400} & \coverentry{14188}{15840} & \coverentry{15253}{16632} & \coverentry{10805}{16800} & \coverentry{1288}{17325} \\
\coverentry{1272}{17600} & \coverentry{5039}{18480} & \coverentry{11887}{18900} & \coverentry{979}{19008} & \coverentry{17429}{19800} \\
\coverentry{15730}{20160} & \coverentry{12921}{20790} & \coverentry{6032}{21600} & \coverentry{17906}{22176} & \coverentry{2771}{23100} \\
\coverentry{3380}{23760} & \coverentry{23285}{24640} & \coverentry{17262}{25200} & \coverentry{23334}{26400} & \coverentry{20944}{27720} \\
\coverentry{22244}{29700} & \coverentry{7550}{30240} & \coverentry{7822}{30800} & \coverentry{5560}{31680} & \coverentry{27480}{33264} \\
\coverentry{4057}{33600} & \coverentry{28859}{34650} & \coverentry{24257}{36960} & \coverentry{27555}{37800} & \coverentry{6931}{39600} \\
\coverentry{7422}{41580} & \coverentry{26748}{43200} & \coverentry{117}{44352} & \coverentry{14953}{46200} & \coverentry{12869}{47520} \\
\coverentry{35}{50400} & \coverentry{43386}{51975} & \coverentry{21190}{52800} & \coverentry{8832}{55440} & \coverentry{28604}{59400} \\
\coverentry{51290}{60480} & \coverentry{2316}{61600} & \coverentry{8436}{66528} & \coverentry{33019}{69300} & \coverentry{65340}{73920} \\
\coverentry{36658}{75600} & \coverentry{33041}{79200} & \coverentry{9471}{83160} & \coverentry{16431}{92400} & \coverentry{18516}{95040} \\
\coverentry{6783}{100800} & \coverentry{91813}{103950} & \coverentry{60518}{110880} & \coverentry{59889}{118800} & \coverentry{3696}{123200} \\
\coverentry{62596}{133056} & \coverentry{87874}{138600} & \coverentry{116975}{151200} & \coverentry{34401}{158400} & \coverentry{90791}{166320} \\
\coverentry{53591}{184800} & \coverentry{147882}{207900} & \coverentry{62786}{221760} & \coverentry{197451}{237600} & \coverentry{86741}{277200} \\
\coverentry{64669}{302400} & \coverentry{173009}{332640} & \coverentry{295599}{369600} & \coverentry{47173}{415800} & \coverentry{225747}{475200} \\
\coverentry{414022}{554400} & \coverentry{25865}{665280} & \coverentry{110190}{831600} & \coverentry{1031743}{1108800} & \coverentry{1097511}{1663200} \\
\coverentry{1956368}{3326400} &  &  &  &  \\
\end{longtable}
\endgroup

\section{Independent Verification}\label{app:verifier}

The supplementary program \texttt{verify.py} requires Python~3.9 or later and no external packages. From the paper directory, run
\begin{lstlisting}[style=pseudo,numbers=none]
python3 anc/supplementary/verify.py main.tex
\end{lstlisting}
It reads the marked certificate blocks, checks distinct moduli, valid residues, minimum modulus and actual LCM, and marks every covered residue in $0,\ldots,L-1$. All four lists pass. Its checks use explicit exceptions and remain active under Python's optimization option.

The supplement also retains reduced certificates and private-point witnesses. A private point is covered by exactly one class; finding one for every retained class certifies irredundancy. These subsidiary checks are independent of the LCM bounds.

\section{Experimental Details and Supplementary Comparisons}\label{app:timing}

\subsection{Execution and budget accounting}\label{app:experiment-details}

The kernels use $32$ perturbation steps and at most $80$ coordinate rounds. Before exceeding a phase allowance, the wrapper interrupts and retains its best state with all active residues assigned, discarding unfinished ruin--recreate proposals. Each insertion phase has a random fallback assignment. At the first covering it stops, assigning arbitrary residues to unused divisors. Thus every output contains the full divisor set. These rules define the new experiments, not the historical executable.

All $180$ assignments and phase/trajectory counts passed independent checks: $75$ coverings and $105$ budget exhaustions. No phase reached the safety cap of $100000$ trajectory starts. Complete start, completion and interruption counts accompany the data.

Runs were serial, with one worker, on an Apple M4 MacBook Pro (24~GB, macOS~27.0.1), using Apple Clang~21.0.0 and \texttt{-std=c++17 -O3 -DNDEBUG -pthread}. CPU and wall timers include setup and initialization; the host was not exclusively reserved. The directory \texttt{supplementary/experiments\_20261031} supplies the frozen protocol, source/executable hashes, commands, assignments and logs.

\subsection{Initial-study attribution and bank diagnostics}\label{app:initial-diagnostics}

At both $L=30240$ and $75600$, five of ten runs per half-insertion configuration succeed during common insertion. Of the other five at $30240$, Adjacent, Global and Bank-1 rescue two each, and Bank-20 rescues three; at $75600$, each rescues one. These conditional counts use the same saved states and are not additional independent samples.

Replays of all $30$ Bank-20 runs exactly reproduce their assignments and non-timing counters. Maximum bank sizes on the three instances are $13,11,9$; nonelite parents are used in $9,4,4$ runs. Replay timings are excluded from Table~\ref{tab:controlled-timing}. In same-random-seed pairs, Bank-20 has lower/equal/higher gaps than Bank-1 in $5/4/1$, $1/8/1$, and $0/10/0$ cases. At $20160$, one Bank-20 run has $72$ extra gaps despite five small paired improvements; the bank medians are $29.5$ versus $21$. There is no consistent larger-bank advantage.

For $0/10$, $3/10$ and $8/10$, the $95\%$ Wilson intervals are approximately $[0,0.278]$, $[0.108,0.603]$ and $[0.490,0.943]$. All row-wise intervals and individual gaps are supplied. The small samples do not establish a general ranking; wall-time matching and projected initialization are untested.

The original bank's exposure audit covers $85$ epochs: $79$ complete one trajectory and interrupt a second, one completes two including success, and five succeed on their first. Every epoch restarts parent selection, so all $165$ starts use positions $1$ or $2$. Schedule cycling also restarts, never reaching the third, mixed schedule. Retained bank size therefore overstates actual parent exposure; the initial capacity comparison is limited to these rules. Section~\ref{sec:followup} tests revised sampling and schedule cycling.

\subsection{Larger-period comparison}\label{app:larger-tests}

The following table gives every larger-period comparison discussed in Section~\ref{sec:followup}. All counts and statistics follow the definitions in Section~\ref{sec:controlled}.

\begin{table}[!htb]
\centering
\caption{Exploratory larger-period tests, $N=5$ per row. Uncovered-residue summaries include all runs; scan medians include only successful runs, with a dash for no success. The scan limit $Q$ is per run. Equal scan counts at different $L$ do not represent equal arithmetic work.}
\label{tab:larger-periods}
\small
\begin{tabular}{@{}lrcrr@{}}
\toprule
Method & \shortstack{Scan limit\\$Q/10^3$} & \shortstack{Successes\\$s/N$} & \shortstack{Uncovered residues\\median [min, max]} & \shortstack{Scans to success\\median ($10^3$)}\\
\midrule
\multicolumn{5}{c}{$m=10,\ L=302400$}\\
Ordered-Whole & 15 & 0/5 & $3741\,[2780,6390]$ & ---\\
Segmented & 15 & 0/5 & $639\,[561,4917]$ & ---\\
Global & 15 & 0/5 & $1880\,[576,4745]$ & ---\\
Ordered-Whole & 60 & 0/5 & $2780\,[1431,3169]$ & ---\\
Segmented & 60 & 0/5 & $775\,[541,1871]$ & ---\\
Global & 60 & 0/5 & $826\,[631,1714]$ & ---\\
\midrule
\multicolumn{5}{c}{$m=11,\ L=3326400$}\\
Ordered-Whole & 15 & 0/5 & $12765\,[215,23453]$ & ---\\
Segmented & 15 & 0/5 & $1014\,[33,4479]$ & ---\\
Global & 15 & 0/5 & $3814\,[330,4352]$ & ---\\
Ordered-Whole & 60 & 1/5 & $712\,[0,1462]$ & 33.9\\
Segmented & 60 & 5/5 & $0\,[0,0]$ & 41.6\\
Global & 60 & 2/5 & $1228\,[0,2295]$ & 26.6\\
\bottomrule
\end{tabular}
\end{table}

\subsection{Supplementary comparison tables}\label{app:order-control}

Table~\ref{tab:ordering} compares initialization orders; Table~\ref{tab:controlled-timing} reports initial-study timing and trajectory counts. A trajectory consists of initialization, coordinate search and perturbation steps (Appendix~\ref{app:trajectory-kernel}); only completed trajectories enter the counts. Successful-run timing is conditional on success, not expected time to success. Follow-up timings and per-parent work are supplied in machine-readable form.

\begin{table}[!htb]
\centering
\caption{Initialization-order control: $N=10$ runs per row, each with scan limit $Q=240000$. Only Ordered-Whole is new; the other rows are initial-study comparators. Uncovered-residue summaries include all runs; scan medians include only successful runs. A dash denotes no success.}
\label{tab:ordering}
\small
\begin{tabular}{@{}lcrr@{}}
\toprule
Method & \shortstack{Successes\\$s/N$} & \shortstack{Uncovered residues\\median [min, max]} & \shortstack{Scans to success\\median ($10^3$)}\\
\midrule
\multicolumn{4}{c}{$m=8,\ L=20160$}\\
Whole & 0/10 & $221.5\,[122,362]$ & ---\\
Ordered-Whole & 0/10 & $151.5\,[96,206]$ & ---\\
Segmented & 0/10 & $60\,[12,102]$ & ---\\
\midrule
\multicolumn{4}{c}{$m=8,\ L=30240$}\\
Whole & 0/10 & $127.5\,[85,302]$ & ---\\
Ordered-Whole & 0/10 & $87.5\,[35,142]$ & ---\\
Segmented & 8/10 & $0\,[0,50]$ & 180.0\\
\midrule
\multicolumn{4}{c}{$m=9,\ L=75600$}\\
Whole & 0/10 & $1096.5\,[633,1881]$ & ---\\
Ordered-Whole & 0/10 & $421.5\,[168,528]$ & ---\\
Segmented & 7/10 & $0\,[0,345]$ & 192.0\\
\bottomrule
\end{tabular}
\end{table}

\begin{table}[H]
\centering
\caption{Initial-study timing: $N=10$ runs per row, each with scan limit $Q=240000$. Uncovered-residue summaries, completed-trajectory medians and process CPU-time medians use all ten runs. Interrupted trajectories are excluded from the completed counts. CPU to success is the median process CPU time at first coverage over successful runs only, with initialization included; a dash denotes no success. All times are in seconds.}
\label{tab:controlled-timing}
\small
\begin{tabular}{@{}lcrrrr@{}}
\toprule
Method & \shortstack{Successes\\$s/N$} & \shortstack{Uncovered residues\\median [min, max]} & \shortstack{Completed\\trajectories} & \shortstack{CPU time\\(s)} & \shortstack{CPU to\\success (s)}\\
\midrule
\multicolumn{6}{c}{$m=8,\ L=20160$}\\
Whole & 0/10 & $221.5\,[122,362]$ & 30 & 1.756 & ---\\
Segmented & 0/10 & $60\,[12,102]$ & 99 & 1.978 & ---\\
Adjacent & 0/10 & $32.5\,[12,82]$ & 65.5 & 1.684 & ---\\
Global & 3/10 & $33\,[0,79]$ & 57 & 1.699 & 1.018\\
Bank-1 & 2/10 & $21\,[0,87]$ & 54 & 1.697 & 1.187\\
Bank-20 & 2/10 & $29.5\,[0,91]$ & 54 & 1.690 & 1.208\\
\midrule
\multicolumn{6}{c}{$m=8,\ L=30240$}\\
Whole & 0/10 & $127.5\,[85,302]$ & 25 & 2.394 & ---\\
Segmented & 8/10 & $0\,[0,50]$ & 74 & 2.164 & 2.108\\
Adjacent & 7/10 & $0\,[0,96]$ & 49 & 1.267 & 1.068\\
Global & 7/10 & $0\,[0,57]$ & 49 & 1.265 & 1.069\\
Bank-1 & 7/10 & $0\,[0,93]$ & 49 & 1.268 & 1.068\\
Bank-20 & 8/10 & $0\,[0,93]$ & 49 & 1.268 & 1.102\\
\midrule
\multicolumn{6}{c}{$m=9,\ L=75600$}\\
Whole & 0/10 & $1096.5\,[633,1881]$ & 20 & 6.203 & ---\\
Segmented & 7/10 & $0\,[0,345]$ & 83 & 6.116 & 6.087\\
Adjacent & 6/10 & $0\,[0,423]$ & 53 & 3.616 & 3.052\\
Global & 6/10 & $0\,[0,382]$ & 51.5 & 3.422 & 3.042\\
Bank-1 & 6/10 & $0\,[0,397]$ & 51.5 & 3.447 & 3.046\\
Bank-20 & 6/10 & $0\,[0,397]$ & 51.5 & 3.420 & 3.042\\
\bottomrule
\end{tabular}
\end{table}

\section*{Acknowledgements}

This work was supported by ShanghaiTech AI Initiative (Grant No.\ AI2026B16). The author thanks Hao Chen, Lian Duan, Boqing Xue, Jiayu Zhai, and Shiliang Zhang for helpful support, discussions, and encouragement.

\bigskip
\noindent\begin{minipage}{\textwidth}
\footnotesize
\noindent J.~Zhang,
\textsc{Institute of Mathematical Sciences\\
ShanghaiTech University\\
No.\ 393 Middle Huaxia Road, Pudong New District,\\
Shanghai, China}\par
\noindent\textsc{Email:} \href{mailto:zhangjh2023@shanghaitech.edu.cn}{\texttt{zhangjh2023@shanghaitech.edu.cn}}
\end{minipage}

\end{document}